\documentclass[11pt]{amsart}
\usepackage{graphicx}
\usepackage{color}\usepackage{amsmath,amsthm,amssymb,amsfonts}
\usepackage[english]{babel}
\usepackage{amsmath}
\usepackage{amssymb}
\newtheorem{teo}{Theorem}[section]
\newtheorem{cor}[teo]{Corollary}

\newtheorem{prop}[teo]{Proposition}
\theoremstyle{definition}

\theoremstyle{remark}

\numberwithin{equation}{section}

\newcommand{\To}{\longrightarrow}

\newcommand{\C}{\mathbb{C}}
\newcommand{\D}{\mathbf{D}}
\newcommand{\B}{\mathcal{B}}
\newcommand{\dem}{\noindent{\textsf{Proof.} }}

\renewcommand{\phi}{\varphi}
\newcommand{\eps}{\varepsilon}
\author[A. Miralles]{Alejandro Miralles}
\address{Alejandro Miralles. Departament de Matem\`atiques and Instituto Universitario de Matem\'{a}ticas y Aplicaciones de Castell\'{o}n (IMAC), Universitat Jaume I de Castell\'{o} (UJI), Castell\'{o}, Spain. \emph{e}.mail:
mirallea@uji.es}

\thanks{The author was partially supported by the projects MTM2014-53241-P, P1-1B2014-35 and AICO/2016/030.}


\begin{document}

\begin{abstract}
The space of Bloch functions on bounded symmetric domains is extended by considering Bloch functions $f$ on the unit ball $B_E$ of finite and infinite dimensional complex Banach spaces in two different ways: by extending the classical Bloch space considering the boundness of $(1-\|x\|^2) \|f'(x)\|$ on $B_E$ and by preserving the invariance of the correspondiing seminorm when we compose with automorphisms $\phi$ of $B_E$. We study the connection between these spaces proving that they are different in general and prove that all bounded analytic functions on $B_{E}$ are Bloch functions in both ways. 
\end{abstract}

\title[Bloch functions on the unit ball of a Banach space]{ Bloch functions on the unit ball of a Banach space}

\maketitle

\section*{Introduction}
The classical Bloch space $\mathcal{B}$ of analytic functions on the open unit disk $\D$ of $\C$ plays an important role in geometric function theory and it has been studied by many authors. K. T. Hahn and R. M. Timoney  extended the notion of Bloch function by considering bounded homogeneous domains in  $\mathbb{C}^n$, such as the unit ball $B_n$ and the polydisk $\D^n$ (see \cite{H,T80,T80b}). O. Blasco, P. Galindo and A. Miralles extended the notion to the infinite dimensional setting by considering Bloch functions on the unit ball of an infinite dimensional Hilbert space (see \cite{BGM, BGLM, BGLM2}) and C. Chu, H. Hamada, T. Honda and G. Kohr considered Bloch functions on bounded symmetric domains which may be also infinite dimensional (see \cite{CHHK}). 

In this article, we will deal with a finite or infinite dimensional complex Banach space $E$ and we will consider two possible extensions of the classical Bloch space. The first one extends the classical Bloch space by considering \textit{the natural Bloch space $\B_{nat}(B_E)$} of holomorphic functions $f$ on $B_E$ such that $\|f\|_{nat} =\sup_{x \in B_E} (1-\|x\|^2)\|f'(x)\| < \infty$. The second one extends the space defined in \cite{CHHK} by considering \textit{the invariant Bloch space $\B_{inv}(B_E)$}  of holomorphic functions $f$ on the unit ball $B_E$ of a complex Banach space $E$ such that $\|f\|_{inv}=\sup_{\phi \in Aut(B_E)} \|(f \circ g)'(0)\| < \infty$. The only known case where $\|\cdot\|_{nat}$ and $\|\cdot\|_{inv}$ are equivalent seminorms and $\B_{nat}(B_E)=\B_{inv}(B_E)$ is when $E$ is a finite or infinite dimensional Hilbert space (see \cite{BGM,T80}). We will prove that there are spaces $E$ satisfying $\B_{inv}(B_E) \subsetneq \B_{nat}(B_E)$ and other ones such that $\B_{nat}(B_E) \subsetneq \B_{inv}(B_E)$. 
Finally we will give a Schwarz-type lemma for complex Banach spaces and will prove that the space of bounded analytic functions on $B_E$ given by $H^{\infty}(B_E)$ is strictly contained in both $\B(B_E)$ and $\B_{nat}(B_E)$.
 

%

\section{Background}

\subsection{The classical Bloch space}

The classical \textit{Bloch space} $\B$ (see \cite{P70}) is the space of analytic functions $f: \D \To \C$ satisfying
$$\|f\|_\B=\sup_{z \in \D} (1-|z|^2) |f'(z)| < \infty$$

endowed with the norm
$$\|f\|_{Bloch}=|f(0)|+\|f\|_\B < \infty$$

\noindent so that $(\B,\| \cdot \|_{Bloch})$ becomes a Banach space. The seminorm $\| \cdot \|_{\B}$ is invariant by automorphisms, that is, $\|f \circ \phi \|_{\B}=\|f\|_{\B}$ for any $f \in \B$ and $\phi \in Aut(\D)$. Recall that 
$$H^{\infty}=\{f:\D \to \C : f \mbox{ is holomorphic and bounded }\}$$
\noindent is a Banach space endowed with the sup-norm $\|f\|_{\infty}=\sup_{z \in \D} |f(z)|$. It is well-known (see for instance \cite{Z07}) that:
\begin{prop} \label{prop hinf y bloch}
 $H^{\infty}$ is properly contained in $\mathcal{B}$ and $\|f\|_{\B} \leq \|f\|_{\infty}$ for any $f \in H^{\infty}$.
\end{prop}
 For further information and references about the classical Bloch space $\mathcal{B}$, the reader is referred  to \cite{ACP, Z07}. 

\subsection{Holomorphic functions on $B_E$ and the pseudohyperbolic distance}

We will denote by $E, F$ complex Banach spaces. Given $x \in E$ and $r>0$, we will denote by $B(x,r)$ the ball given by $y \in E$ such that $\|y-x\| < r$. We will denote by $B_E$ the open unit ball $B(0,1)$ of $E$. A function $f:B_{E} \to F$
is said to be holomorphic if it is Fr{\'e}chet differentiable at every  $x \in B_{E}$ or, equivalently, if $f(x)=\sum_{n=1}^{\infty} P_n(x)$ for all $x \in B_{E}$, where $P_n$ is an $n-$homogeneous polynomial, that is, the restriction to the diagonal of a continuous $n-$linear form on the $n$-fold space $E\times\dots\times E$ into $F$. We will denote by $H(B_E,F)$ the space of holomorphic functions from $B_E$ into $F$. If $F=\C$, we just denote the space by $H(B_E)$. For further information on holomorphic functions on complex Banach spaces, see \cite{D} or \cite{M}. \smallskip 

The space $H^{\infty}(B_{E})$ is given by $\{f:B_{E} \to \mathbb{C}: f\hbox{ \textit{is holomorphic
and bounded }} \}$ and it becomes a Banach space when endowed with
the sup-norm $\|f\|_{\infty}=\sup\{|f(x)|:x\in B_{E}\}$. 

\subsection{The pseudohyperbolic and the hiperbolic distance on $B_E$} \label{pseud}
The pseudohyperbolic distance $\rho$ for $z,w \in \D$ is given by
$$\rho(z,w)=\left| \frac{z-w}{1-\bar{z}w} \right|.$$

\noindent The pseudohyperbolic distance $\rho$ for
$x, y \in B_E$ is given by 
$$\rho(x,y)=\sup\{ |f(x)| : \ f \in H^{\infty}(B_E), \ \|f\| \leq 1, \ f(y)=0 \}.$$

\noindent We recall the following well-known  results:
\begin{prop}
Let $E$ be a complex Banach space and $x,y \in B_E$. Then:
\begin{itemize}
\item[a)] $\rho(f(x),f(y)) \leq \rho(x,y)$ for any $f \in H^{\infty}(B_E)$ such that $\|f\|_{\infty} \leq 1$.
\item[b)] $\rho(\phi(x),\phi(y)) \leq \rho(x,y)$ for all holomorphic mappings $\phi: B_E \to B_E$. The equality is satisfied if and only if $\phi \in Aut(B_E)$.
\end{itemize}
\end{prop}

The hyperbolic distance $\beta$ for $x,y \in B_E$ is given by
$$\beta(x,y)=\frac{1}{2} \log \left( \frac{1+\rho(x,y)}{1-\rho(x,y)} \right).$$

For bounded symmetric domains $B_E$, it was proved that any Bloch function $f$ on $B_E$ is Lipschitz for the hyperbolic distance (see \cite{BGLM} and \cite{CHHK}), that is, there exists $M >0$ such that for any $x,y \in B_E$,
$$|f(x)-f(y)| \leq M \beta(x,y).$$

\subsection{The space of Bloch functions on bounded symmetric domains} \label{defis}

The study of Bloch functions on bounded symmetric domains of $\C^n$ was extended by Hahn \cite{H} and Timoney by using the Bergman metric (see \cite{T80,T80b}). In particular, this study includes the unit euclidean ball $B_n$ and the polydisc $\D^n$. The study of Bloch functions on bounded symmetric domains of infinite dimensional Banach spaces was introduced by Blasco, Galindo and Miralles (see \cite{BGM}) for the Hilbert case and by Chu, Hamada, Honda and Kohr for general bounded symmetric domains by means of the Kobayashi metric (see \cite{CHHK}). If we consider these domains as the unit ball $B_E$ of a $JB^*-$triple $E$, the corresponding Bloch space is the set of holomorphic functions on $B_E$ which satisfy that $\sup_{x \in B_E} Q_f(z) < \infty$. In this case, 
$$Q_f(z)=\sup_{x \in X \setminus \{0\} } \frac{|f'(z)(x)|}{k(z,x)}$$
\noindent and 
$k(z,x)$ is the infinitesimal Kobayashi metric for $B_E$ (see \cite{CHHK} for more details). It was proved that
$$\sup_{x \in B_E} Q_f(x) = \sup_{\phi \in Aut(B_E)} \|(f \circ \phi)'(0)\|$$
so the Bloch space $\B(B_E)$ on a bounded symmetric domain can be described in terms of the automorphisms of $B_E$. The authors also proved that $H^{\infty}(B_E) \subsetneq \B(B_E)$.

\subsection{The automorphisms on $B_E$}

We will denote by $Aut(B_E)$ all the automorphisms of $B_E$, that is, all the bijective biholomorphic maps $\phi: B_E \to B_E$. It is well-known that if $B_E$ is a bounded symmetric domain (including the unit ball of a Hilbert space and the finite or infinite dimensional polydisc) then they are homogeneous, that is, they act transitively on $B_E$. Hence, if $B_E$ is a bounded symmetric domain, then $\{ \phi(0) : \phi \in \mbox{Aut}(B_E) \} = B_E$. Kaup and Upmeier (see \cite{KU}) proved that:
\begin{prop} \label{KU}
If $E$ is a complex Banach space and $B_E$ is its open unit ball, then the set $V=\{\phi(0) : \phi \in Aut(B_E) \}$ is a closed subspace of $E$ and $B_E \cap V$ is a bounded symmetric domain in $V$.
\end{prop}

Hence, it is clear that
\begin{prop}
Let $E$ be a complex Banach space and $B_E$ its open unit ball and consider $V=\{\phi(0) : \phi \in Aut(B_E) \}$. Then $B_V$ is a bounded symmetric domain and $Aut(B_V)=\{\phi|_{B_E \cap V} : \phi \in Aut(B_E)\}$.
\end{prop}

\dem Notice that for any $\phi \in Aut(B_E)$ we have that
$\phi(B_E \cap V) =B_E \cap V$ since for any $x \in B_E \cap V$ we have that $x=\psi(0)$ for some $\psi \in Aut(B_E)$. Hence, $\phi(x)=(\phi \circ \psi)(0)$ which clearly belongs to $B_E \cap V$ because $\phi \circ \psi \in Aut(B_E)$. It is clear that $\phi|_{B_E \cap V}$ is also surjective since for any $y \in B_E \cap V$ there exists $\psi_1 \in Aut(B_E)$ such that $\psi_1(0)=y$. Since $\phi^{-1} \circ \psi_1 \in Aut(B_E)$, take $x = \phi^{-1}(\psi_1(0)) \in B_E \cap V$ and it is clear that $\phi(x)=y$ so we are done. \qed \bigskip

\section{The space of Bloch functions on $B_E$}

\subsection{Two different definitions}

Let $E$ be a complex Banach space and consider its open unit ball denoted by $B_E$. Bearing in mind the definition of the classical Bloch space taking the supremum of $(1-|z|^2) |f'(z)|$ for $z \in \C, |z| <1$, we can extend it for $f \in H(B_E)$ by defining what we call \textit{the natural Bloch seminorm}
$$\|f\|_{nat}=\sup_{x \in B_E} (1-\|x\|^2)\|f'(x) \|$$

\noindent where $f'(x) \in E^*$ denotes the derivative of $f$ at the point $x$. The space $\B_{nat}(B_E)$ is given by
$$\B_{nat}(B_E)=\{f \in H(B_E) : \|f\|_{nat} < \infty \}.$$
\noindent It is clear that $\| \cdot \|_{nat}$ is a seminorm for $\B_{nat}(B_E)$ and this space can be endowed with the norm $\|f\|_{nat-Bloch}=|f(0)|+\|f\|_{nat}$. It is easy that $(\B_{nat}(B_E),\| \cdot \|_{nat-Bloch})$ is a Banach space. \smallskip

On the other hand, bearing in mind the definition of the Bloch space in \cite{CHHK} and to preserve the invariance of the corresponding seminorm when composing with an automorphism, we define for $f \in H(B_E)$, the Bloch semi-norm by
$$\|f\|_{inv}=\sup_{\phi \in Aut(B_E)} \|(f \circ \phi)'(0) \|$$

\noindent and the space $B_{inv}(B_E)$ will be given by
$$B_{inv}(B_E)=\{f \in H(B_E) : \|f\|_{inv} < \infty \}.$$
\noindent It is clear that  $\| \cdot \|_{inv}$ fails to be a norm when we add up $|f(0)|$ if we deal with $B_E$ for general Banach spaces $E$ since Proposition \ref{KU} does not assure that $V=B_E$ if $B_E$ is not a bounded symmetric domain. So we consider the quotient space
$$\B_{inv}(B_E) = B_{inv}(B_E) \slash \sim$$ 
\noindent where $f \sim g$ if and only if $\|f\|_{inv}= \|g\|_{inv}$. We endow this space with the norm $\|f\|_{inv-Bloch}=\|f\|_{inv}$ and it becomes a Banach space. If we deal with bounded symmetric domains $B_E$, the corresponding space of Bloch functions with the invariant seminorm coincide with the one defined in \cite{CHHK} and the norms are equal up to the constant $|f(0)|$.

As we have mentioned, these seminorms are equivalent if $E$ is a finite or infinite dimensional Hilbert space, so we have that $\B_{nat}(B_E)=\B_{inv}(B_E)$ in this case. 
In the case of $B_E$ which are bounded symmetric domains, it was proved (see Corollary 3.5 in \cite{CHHK}):
\begin{cor}
If $B_E$ is a bounded symmetric domain, then for any $x \in B_E$ we have
$$\|f'(x)\| \leq \frac{\|f\|_{B}}{1-\|x\|^2}.$$
\end{cor}

Hence, we have that
\begin{prop} \label{cont_sym}
Let $E$ be a complex Banach space such that $B_E$ is a bounded symmetric domain. Then, 
$\B_{inv}(B_E) \subseteq \B_{nat}(B_E)$ and $\|f\|_{nat} \leq \|f\|_{inv}$.
\end{prop}

In this section, we will give examples where these spaces are different even for some bounded symmetric domains. Indeed, for general Banach spaces $E$ we will show that it is not true that $\B_{inv}(B_E) \subseteq \B_{nat}(B_E)$ or $\B_{nat}(B_E) \subseteq \B_{inv}(B_E)$.

\subsection{The case $E=(\C^n,\| \cdot \|_{\infty})$ and $E=c_0$}
Let $E=(\C^n,\| \cdot \|_{\infty})$ or $E=c_0$, whose open unit ball is the so-called (finite or infinite dimensional) polydisc, which is usually denoted by $\D^n$ and $B_{c_0}$ respectively. For any $f \in H(B_E)$ and $x \in B_E$ we have that $f'(x)$ belongs to $\ell_1^n$ or $\ell_1$ respectively, so we can identify $f'(x)$ by $\left( \frac{\partial f}{\partial x_1}(x), \frac{\partial f}{\partial x_2}(x), \cdots \right)$ and $\|f'(x)\|=\sum_{k=1}^{n} \left| \frac{\partial f}{\partial x_k}(x) \right|$, where $n$ can be finite or infinite.

\subsubsection{Automorphisms on $E=(\C^n,\| \cdot \|_{\infty})$ and $E=c_0$}

Bloch functions on the finite or infinite polydisc were studied in \cite{T80} and \cite{CHHK} respectively. In these works, the authors deal with the natural Bloch seminorm $\|\cdot\|_{nat}$ and compare it with the invariant Bloch seminorm $\| \cdot\|_{inv}$ in the case of a Hilbert space $E$. Now we prove that the spaces defined with each of these seminorms are different even if $B_E$ is the bidisc $\D^2$. \medskip


Consider for any $z \in \D$ the automorphism $\phi_z: \D \to \C$ given by
$$\phi_z(w)=\frac{w-z}{1-\bar{z}w},$$
\noindent which satisfies that $\phi_z'(0)=-(1-|z|^2)$. It is well-known that any $\phi \in Aut(\D)$ is given by $\phi(w)=e^{i \alpha} \phi_z(w)$ for some $z \in \D$ and $\alpha \in [0,2 \pi[$. \medskip

Now let $x=(x_1,x_2,\cdots) \in B_E$  and consider the automorphism $\phi_x: B_{E} \to B_{E}$ given by 
$$\phi_x(y)=\left(\frac{x_1-y_1}{1-\bar{x_1}y_1}, \frac{x_2-y_2}{1-\bar{x_2}y_2},\cdots\right) \cdot$$
\noindent It is well-known that any $\phi \in Aut(B_E)$ is given by $\phi=(\phi_{1},\phi_{2}, \cdots)$ where $\phi_k \in Aut(\D)$ (see \cite{R} and \cite{FJ} for the finite and infinite dimensional case respectively).

We are interested in the calculation of $\|(f \circ \phi)'(0)\|$ for $\|f\|_{inv}$, where $f: B_E \to \C$ is a holomorphic function and $\phi$ is an automorphism of $B_E$, so we can consider, without loss of generality, that $\phi=\phi_x$ since $|e^{i\alpha} z|=|z|$ for $z \in \D$, $\alpha \in [0,2 \pi[$ and $(e^{i \alpha_1} z_1,e^{i \alpha_2} z_2, \cdots) \in B_E$ for any $(z_1,z_2,\cdots) \in B_E$, $\alpha_1,\alpha_2,\cdots \in [0,2\pi[$.

\begin{prop}
Let $f \in H(B_E)$ where $E=(\C^n,\| \cdot \|_{\infty})$ or $E=c_0$. Then,
$$\|(f \circ \phi_x)'(0)\|=\sum_{k=1}^n (1-|x_k|^2)\left| \frac{\partial f}{\partial x_k}(x)\right|,$$
\noindent where $n$ can be finite or infinite.
\end{prop}

\dem It is clear that 
$$\phi_x'(0)= \left( \begin{array}{ccc}
-(1-|x_1|^2) & 0 & \cdots \\
0 & -(1-|x_2|^2) & \cdots \\
0 & 0 & \cdots  \\ 
\cdots & \cdots & \cdots \end{array} \right)$$

\noindent Since $f \circ \phi_x$ is well-defined on $B_{E}$ and bearing in mind that $\phi_x(0)=x$, we have 
$$(f \circ \phi_x)'(0)=f'(\phi_x(0)) \circ \phi_x'(0)= \left( \begin{array}{ccc}
\frac{\partial f}{\partial x_1}(x) & \frac{\partial f}{\partial x_2}(x)& \cdots  \end{array} \right) \circ \left( \begin{array}{ccc}
-(1-|x_1|^2) & 0 & \cdots \\
0 & -(1-|x_2|^2) & \cdots \\
0 & 0 & \cdots \end{array} \right)$$
\noindent so
$$\|(f \circ \phi_x)'(0)\|=\sum_{k=1}^n (1-|x_k|^2)\left| \frac{\partial f}{\partial x_k}(x)\right|,$$
\noindent where $n$ can be finite or infinite. \qed \bigskip

Hence
\begin{cor}
For any $f \in \B_{inv}(B_E)$, we have that
$$\|f\|_{inv}=\sup_{x \in B_{E}} \|(f \circ \phi_x)'(0) \|=\sup_{x \in B_E} \sum_{k=1}^n (1-|x_k|^2) \left|\frac{\partial f}{\partial x_k}(x) \right|$$
\noindent where $n$ can be finite or infinite.
\end{cor}


%
%
%
%

\begin{prop} \label{countex1}
Let $E=(\C^2,\| \cdot \|_{\infty})$ and $B_E=\D^2$ the bidisc. Then $\B_{inv}(\D^2)) \subsetneq \B_{nat}(\D^2)$.
\end{prop}

\dem It is clear that $\B_{inv}(\D^2) \subseteq \B_{nat}(\D^2)$ by Proposition \ref{cont_sym}. To prove that these spaces are different, consider $f(z,w)=(w+1) \log (z-1)$. Then,
$\frac{\partial f}{\partial z}=\frac{w+1}{z-1}$ and $\frac{\partial f}{\partial w}=\log (z-1)$. Notice that
$$\|f\|_{nat}=\sup_{(z,w) \in \D^2} (1-\sup\{|z|,|w|\}^2) \left( \left| \frac{\partial f}{\partial z}(z,w) \right| + \left| \frac{\partial f}{\partial w}(z,w) \right| \right)$$
\noindent so
$$\|f\|_{nat} \leq \sup_{|z| < 1} (1-|z|^2) \left( \left|\frac{w+1}{z-1} \right| + \left| \log(z-1) \right| \right) \leq 2 \sup_{|z| < 1} (1-|z|) \left( \frac{2}{1-|z|} + \left| \log(z-1) \right| \right) \leq$$
$$4 +  \sup_{|z| < 1} |z-1| |\log(z-1)|$$

It is clear that $w \log w$ is bounded on the set of complex numbers $w$ such that $|w| \leq 2$ since
$|w \log w| \leq |w| |\log|w|+i \arg w| \leq |w|( \log|w|+2 \pi)$ and $t \log t \rightarrow 0$ when $t \rightarrow 0$ so there exists a constant $M >0$ such that $\|f\|_{nat} \leq 4+M < \infty$ and concude that $f \in \B_{nat}(\D^2)$. \smallskip

However, evaluating in $w=0$,
$$\|f\|_{inv} =\sup_{(z,w) \in \D^2} (1-|z|^2) \left| \frac{\partial f}{\partial z}(z,w)\right| +(1-|w|^2) \left| \frac{\partial f}{\partial w}(z,w) \right| \geq$$
$$\sup_{(z,w) \in \D^2} (1-|z|) \left| \frac{w+1}{z-1}  \right| +(1-|w|) \left| \log(z-1) \right| \geq \sup_{|z| <1} (1-|z|) \left| \frac{1}{z-1}  \right| + \left| \log(z-1) \right|$$
\noindent and taking $z_n=1-1/n$, since  
$$(1-|z_n|) \left| \frac{1}{z_n-1}  \right| + \left| \log(z_n-1) \right|=\frac{1}{n} \frac{1}{\frac{1}{n}}+|\log \left( \frac{-1}{n} \right)|\geq 1+\left|\log \left( \frac{1}{n}\right)\right| \rightarrow \infty$$
\noindent when $n \rightarrow \infty$, so $f \notin \B_{inv}(\D^2)$. \qed \bigskip

Hence, we have
\begin{cor}
Let $E$ be $(\C^n,\| \cdot\|_{\infty})$ for $n \geq 2$ or $c_0$. Then, $\B_{inv}(B_E) \subsetneq \B_{nat}(B_E)$.
\end{cor}

\dem It is clear that $\B_{inv}(B_E) \subseteq \B_{nat}(B_E)$ because of Proposition \ref{cont_sym}. Consider $f \in H(\D^2)$ given in Proposition \ref{countex1}. The function $g(x_1,x_2,x_3,\cdots)=f(x_1,x_2)$ belongs to $\B_{nat}(B_E)$ but $g \notin \B_{inv}(B_E)$. \qed \bigskip

\subsection{The case $E=L_p(\Omega,\mu)$} \label{sec_lp}

L. L. Stach\'o and E. Vesentini (see \cite{S} and \cite{V}) proved that for measure spaces $E=L_p(\Omega,\mu)$, $1 \leq p < \infty$, $p \neq 2$ and $\mu(\Omega) < \infty$, we have
$$Aut(B_E)=\{ U|_{B_E} : U \mbox { is a surjective linear isometry  of } E \}.$$

Hence, $\phi'(0)=\phi$ and $\phi(0)=0$ for any $\phi \in Aut(B_E)$. \smallskip

We will prove that the behaviour of the unit ball $B_E$ of these spaces is completely different to bounded symmetric domains when we deal with the spaces of Bloch functions on $B_E$.

\begin{prop} \label{counter2}
Let $E=L_p(\Omega,\mu)$, $1 \leq p < \infty$, $p \neq 2$ and $\mu(\Omega) < \infty$.  Then $\B_{nat}(B_E) \subsetneq \B_{inv}(B_E)$.
\end{prop}

\dem Let $\phi \in Aut(B_E)$. Since $\phi$ is the restriction of a surjective linear isometry to $B_E$, we have
$$\|(f \circ \phi)'(0)\|=\| f('(\phi(0)) \circ \phi'(0)\|= \|=\|f'(0) \circ \phi'(0) \|=\|f'(0)\| \leq \sup_{x \in B_E}(1-\|x\|^2) \|f'(x)\|$$
\noindent so $\B_{nat}(B_E) \subseteq \B_{inv}(B_E)$. \smallskip

However, J. M. Ansemil, R. Aron and S. Ponte (see \cite{A} and \cite{AAP}) proved that given any two disjoint balls in an infinite dimensional complex Banach space $E$,
there exists an entire function on $E$ which is bounded on one and unbounded on the other. We consider the balls $B_1=\frac{1}{2} B_E=\{ x \in E : \|x\| < \frac{1}{2} \}$ and $B_2=B(x_0,\frac{1}{5}):=\{ x \in E : \|x-x_0\| < \frac{1}{5} \}$ for a fixed $x_0 \in E$ such that $\|x_0\|=\frac{3}{4}$. Then, there exists an entire function $f$ on $E$ such that $f|{B_1}$ is bounded and $f|_{B_2}$ is unbounded, so there exists $x_n \subset B_2$ such that $|f(x_n)| \rightarrow \infty$ when $n \rightarrow \infty$. By the Mean Value Theorem (see Theorem 13.8 in \cite{M}) we have that
$$|f(x_n)-f(0)| \leq \|x_n\| \sup_{0 \leq \lambda \leq 1} \|f'(\lambda x_n)\|,$$
\noindent so 
$$\frac{|f(x_n)| -|f(0)|}{\|x_n\|} \leq \frac{|f(x_n)-f(0)|}{\|x_n\|} \leq \sup_{0 \leq \lambda \leq 1} \|f'(\lambda x_n)\|$$
\noindent and since $\|x_n\|  \leq \frac{3}{4}$, we have that
$$\frac{|f(x_n)|-|f(0)|}{\|x_n\|} \geq \frac{4}{3} \left( |f(x_n)|-|f(0)| \right) \rightarrow \infty$$ \noindent when $n \rightarrow \infty$, so we can take a sequence $(\lambda_n) \subset[0,1]$ such that $ \|f'(\lambda_n x_n)\| \rightarrow \infty$ when $n \rightarrow \infty$. Hence,
$$(1-\| \lambda_n x_n \|^2) \|f'(\lambda x_n)\| \geq (1-\| \lambda_n x_n\|) \|f'(\lambda_n x_n)\| \geq (1-\|x_n\|) \|f'(\lambda_n x_n)\| \geq \frac{1}{4} \|f'(\lambda_n x_n)\|  \rightarrow \infty$$
\noindent when $n \rightarrow \infty$.

Since $B_1,B_2 \subset B_E$, we conclude that $\|f\|_{nat} =\infty$ but, as we have observed above, $\|f\|_{inv} = \|f'(0)\| < \infty$ so  $\B_{nat}(B_E) \subsetneq \B_{inv}(B_E)$. \qed \bigskip

\subsection{Bloch functions and Lipschitz functions for the hyperbolic metric}

As we have mentioned in Subsection \ref{pseud}, any $f \in \B_{inv}(B_E)$ is Lipschitz for the corresponding hyperbolic metric $\beta$ on $B_E$ for any Banach space $E$ such that $B_E$ is a bounded symmetric domain (see \cite{CHHK}).

We will prove that if we deal with $B_E$ which are not bounded symmetric domains, this is no longer true. Consider the spaces $L_p$ considered in subsection \ref{sec_lp}. Then,
\begin{prop}
Let $E=L_p(\Omega,\mu)$ be as above. Then there exists $f \in \B_{inv}(B_E)$ which is not Lipschitz for the corresponding hyperbolic distance $\beta$ on $B_E$.
\end{prop}

\dem Look at the proof of Proposition \ref{counter2}. Take the balls $B_1$, $B_2$, $f$ the function which is defined there and the sequence $(x_n) \subset B_E$ such that $|f(x_ n)| \rightarrow \infty$ when $n \rightarrow \infty$. So $|f(x_n)-f(0)| \geq |f(x_n)|-|f(0| \rightarrow \infty$ when $n \rightarrow \infty$ but 
$$\beta(x_n,0)=\frac{1}{2} \log \left( \frac{1+\|x_n\|}{1-\|x_n\|} \right)$$
\noindent which is bounded on $(x_n)$ since $1-\|x_n\| \geq 1-\frac{3}{4}=\frac{1}{4}$. \qed \bigskip

\section{Bounded functions on $B_E$ are Bloch functions}

In \cite{BGM} and \cite{CHHK} it was proved that $H^{\infty}(B_E) \subsetneq \B(B_E)$ when $E$ is a Hilbert space or $B_E$ is a bounded symmetric domain respectively. In this section, we will prove that this result remains true if we deal with any complex Banach space $E$ and any Bloch space, that is, and $H^{\infty}(B_E) \subsetneq \B_{nat}(B_E)$ and $H^{\infty}(B_E) \subsetneq \B_{inv}(B_E)$.  \smallskip

First we recall the following result which is an application of the Schwarz lemma (see page 641 in [CCG89]).
\begin{teo} \label{schwarz}
If $g:B(x_0,r)\to \mathbb{C}$ is analytic and $ g(x_0)=0,$ then $$|g(y)|\leq ||g||\frac{\|y-x_0\|}{r} \text{ if } \|y-x_0\|< r.$$
\end{teo}

As consequence, it is clear by the definition of the pseudohyperbolic distance that
\begin{cor} \label{ineq1}
Let $E$ be a complex Banach space and $x,y \in B_E$. Then,
\begin{eqnarray} 
\rho(y,x) \leq \frac{\|y-x\|}{r} \mbox{ for all } y \in B(x,r).
\end{eqnarray}
\end{cor}

\begin{prop} \label{propo}
Let $f \in H^{\infty}(B_{E})$ such that $\|f\|_{\infty} \leq 1$. Then, for any $x_0 \in B_{E}$, we have that
$$(1-\|x_0\|) \| f'(x_0) \| \leq 1-|f(x_0)|^2.$$
\end{prop}

\dem First we consider that $\|f\|_{\infty} <1$ and let $x_0 \in B_E$. Applying Corollary \ref{ineq1}, for $r >0$ such that $\|x_0\|+r<1,$ we have that $$\rho(y,x_0)\leq \frac{\|y-x_0\|}{r}~~~~\;\mbox{ for all } y \in B(x_0,r).$$ 

Taking limits when $r \rightarrow (1-\|x_0\|)^{-}$ we get $$\rho(y,x_0)\leq \frac{\|y-x_0\|}{1-\|x_0\|} \text { for all } y\in B(x_0,1-\|x_0\|).$$

Notice that for any $x_0 \in B_{E}$,  $f'(x_0)$ is the functional on $E$ satisfying
$$\lim_{x \rightarrow x_0} \frac{f(x)-f(x_0) -f'(x_0)(x-x_0)}{\|x-x_0\|} =0,$$

\noindent so given $\eps > 0$, there exists $\delta >0$ such that
$$\left| \frac{f(x)-f(x_0) -f'(x_0)(x-x_0)}{\|x-x_0\|} \right| < \eps,$$

\noindent if $\|x-x_0\| \leq \delta$. Without loss of generality, we can choose $\delta$ such that $x_0 + B(x_0,\delta) \subset B_{E}$. For $x$ such that $\|x-x_0\| \leq \delta$ we have that 
$$\left| \frac{f'(x_0)(x-x_0)}{\|x-x_0\|} \right| \leq \eps +  \frac{\left|f(x)-f(x_0)\right|}{\|x-x_0\|} .$$

Choose a sequence $(\eps_n)$ of positive numbers such that $\eps_n \rightarrow 0$ and consider their corresponding $(\delta_n)$. Since
$\|f'(x_0)\| = \sup_{y \in B_{E}} |f'(x_0) (y)|$, we choose vectors $(y_n) \subset B_{E}$, $\|y_n\| \rightarrow 1$ such that $\|f'(x_0)\| = \lim_{n \rightarrow \infty} |f'(x_0) (y_n)|$ and define
$x_n \in B_{E}$ by
$$x_n=x_0+\delta_n \frac{y_n}{\|y_n\|}.$$
\noindent It is clear that $x_n \in B_E$ since $x_0 + B(x_0,\delta) \subset B_{E}$ and
$$\left\| \frac{x_n-x_0}{\|x_n-x_0\|} -y_n \right\|=\left\| \frac{y_n}{\|y_n\|} - y_n \right\| \rightarrow 0$$ 

\noindent when $n \rightarrow \infty$ since $\|y_n\| \rightarrow 1$. Hence,
 
$$\|f'(x_0)\| = \lim_{n \rightarrow \infty} |f'(x_0)(y_n)|= \lim_{n \rightarrow \infty}  \frac{\left| f'(x_0)(x_n-x_0)\right|}{\|x_n-x_0\|}.$$ 

Notice that
$$\frac{\left| f'(x_0)(x_n-x_0)\right|}{\|x_n-x_0\|} \leq \eps_n +  \frac{\left| f(x_n)-f(x_0) \right|}{\|x_n-x_0\|} \leq \eps_n + \left| \frac{f(x_n)-f(x_0)}{1- \overline{f(x_0)} f(x_n)} \right|  \frac{\left| 1-\overline{f(x_0)} f(x_n) \right|}{\|x_n-x_0\|}.$$

Since the pseudohyperbolic distance is contractive for $f$, we have that
$$\frac{\left| f'(x_0)(x_n-x_0)\right|}{\|x_n-x_0\|}  \leq \eps_n + \rho(x_n,x_0)\frac{ \left| 1-\overline{f(x_0)} f(x_n)\right|}{\|x_n-x_0\|}.$$


Next observe that there is no loss of generality in assuming $\delta_n \to 0,$ so $\delta_n \leq 1-\|x_0\|.$ Hence $$\rho(x_n,x_0)\leq \frac{\|x_n-x_0\|}{1-\|x_0\|}.$$
\noindent so
$$\left| \frac{f'(x_0)(x_n-x_0)}{\|x_n-x_0\|} \right| \leq \eps_n + \frac{\|x_n-x_0\|}{1-\|x_0\|} \frac{|1-\overline{f(x_0)} f(x_n)|}{\|x_n-x_0\|}= \eps_n + \frac{|1-\overline{f(x_0)}f(x_n)|}{1-\|x_0\|}.$$

Taking limits when $n \rightarrow \infty$, we have
$$\|f'(x_0)\| \leq \frac{1}{1-\|x_0\|} (1-|f(x_0)|^2)$$

\noindent and we are done. Suppose now that $\|f\|_{\infty}=1$. Then, there exists a sequence of functions $(f_n) \subset H^{\infty}(B_{E})$ (for instance, $f_n(x)=(1-1/n)f(x)$) such that $f_n$ converges uniformly to $f$ on $B_{E}$. We apply the inequality to the functions $(f_n)$ and taking limits when $n \rightarrow \infty$, we are done. \qed

\begin{prop} \label{propo2}
Let $E$ be a complex Banach space and $B_E$ its open unit ball. Then $H^{\infty}(B_E) \subset \B_{nat}(B_E)$ and the map $Id: (H^{\infty}(B_E),\|\cdot\|_\infty) \to (\B_{nat}(B_E),\|\cdot\|_{nat-Bloch})$ is continuous.
\end{prop}

\dem Let $f \in H^{\infty}(B_E) $. Then $f/\|f\|_{\infty}$ has sup-norm $1$ and we can apply Proposition \ref{propo}. Then, for any $x \in B_E$ we have
$$\frac{\|f'(x)\|}{\|f\|_{\infty}} \leq \frac{1}{1-\|x\|} \left( 1- \frac{|f(x)|^2}{\|f\|_{\infty}^2}\right) \leq \frac{1}{1-\|x\|},$$
\noindent so $(1-\|x\|^2) \|f'(x)\|=(1+\|x\|)(1-\|x\|) \|f'(x)\| \leq 2 \|f\|_{\infty}$ and we obtain that $H^{\infty}(B_E) \subset \B_{nat}(B_E)$. Adding up $|f(0)|$ to the left term we obtain that $\|Id(f)\|_{nat-Bloch} \leq 3 \|f\|_{\infty}$, so bounded functions are Bloch functions and the inclusion is continuous. \qed \bigskip

Now we prove that the same result remains true if we deal with $\B_{inv}(B_E)$ instead of $\B_{nat}(B_E)$.
\begin{prop}
Let $E$ be a complex Banach space and $B_E$ its open unit ball. Then $H^{\infty}(B_E) \subset \B_{inv}(B_E)$ and the map $Id: (H^{\infty}(B_E),\|\cdot\|_\infty) \to (\B_{inv}(B_E),\|\cdot\|_{Bloch})$ is continuous.
\end{prop}

\dem Let $f \in H^{\infty}(B_E)$. For any $\phi \in Aut(B_E)$ it is clear that $f \circ \phi \in H^{\infty}(B_E)$ and $\|f \circ \phi\|_{\infty} \leq \|f\|_{\infty}$ since $\phi(B_E) \subset B_E$. So by the proof of Proposition \ref{propo2} we have
$\|(f \circ \phi)'(0)\|  \leq 2\|f\|_{\infty}$ and hence $\|f\|_{inv} \leq 2\|f\|_{\infty}$. We conclude that $H^{\infty}(B_E) \subseteq \B_{inv}(B_E)$ and in addition,

$\|f\|_{inv-Bloch} =|f(0)|+\|f\|_{inv} \leq 3\|f\|_{\infty}$
\noindent so $Id$ is also continuous. \qed \bigskip
Finally we prove that $H^{\infty}(B_E)$ is strictly contained in $\B_{nat}(B_E)$ or $\B_{inv}(B_E)$. 
\begin{prop}
For any inifinite dimensional complex Banach space $E$, we have that $H^{\infty}(B_E) \subsetneq \B_{nat}(B_E)$.
\end{prop}

\dem Let $x_0 \in E$, $\|x_0\|=1$ and let $L \in E^*$ such that $\|L\|=1$ and $L(x_0)=1$. The function $f(x)=\log(1-L(x))$ satisfies that $f \in \B_{nat}(B_E) \setminus H^{\infty}(B_E)$ since 
$$(1-\|x\|^2)\|f'(x)\|=(1-\|x\|^2)\frac{\|L\|}{|1-L(x)|} \leq (1-\|x\|^2)\frac{\|L\|}{1-\|x\|} \leq 2 \|L\|$$
\noindent but there exists $(x_n) \subset B_E$ such that $x_n \rightarrow x_0$ and $\lim_{n \rightarrow }|f(x_n)|=|\log(1-L(x_n))|=\infty$, so $f \notin H^{\infty}(B_E)$. \qed \bigskip

\begin{prop}
For any inifinite dimensional complex Banach space $E$, we have that $H^{\infty}(B_E) \subsetneq \B_{inv}(B_E)$.
\end{prop}

\dem By Proposition \ref{KU} we know that $V=\{\phi(0) : \phi \in Aut(B_E) \}$ is a closed subspace of $E$. If $V=\{ 0 \}$ we are done since any automorphism $\phi$ satisfies $\phi(0)=0$ and then it is the restriction of a linear isometry of $E$ (see Proposition 1 in \cite{D2}). Then for any $f \in H(B_E)$ we have that $\|f\|_{inv}=\|f'(0)\| < \infty$ so any $f \in H(B_E)$ belongs to $\B_{inv}(B_E)$ but it is well-known that there are unbounded holomorphic functions on $B_E$. If $V \neq \{0\}$ there exists a linear map $L: V \to \C$ and $x \in V$ such that $\|L\|=1$ and $\|x\|=1$, $L(x)=\|L\|$. By the Hahn-Banach Theorem, there exists a linear map $L_1$ on $E$ such that $L_1|V=L$ and $\|L_1\|=\|L\|$. We consider the map $f(x)=\log(1-L(x))$. It is clear that $f$ is unbounded on $B_E$ since $B_V \subset B_E$. Now we prove that $f \in \B_{inv}(B_E)$. Notice that for any $\phi \in Aut(B_E)$ we have that $L \circ \phi \in H^{\infty}(B_E)$ since $\phi(B_E) \subset B_E$ and $\|L \circ \phi \|_{\infty} \leq \|L\|$=1. Let $h(z)=\log(1-z)$, so $f(x)=(h \circ L)(x)$. For any $\phi \in Aut(B_E)$ we have that
$$\| (f \circ \phi)'(0) \| =\|(h \circ L \circ \phi)'(0)\| \leq \|h'(L ( \phi (0))\| \|(L \circ \phi)'(0) \| \leq \frac{\|h\|_{B}}{1-|L(\phi(0))|^2} \|(L \circ \phi)'(0) \| $$ 
\noindent where $\|h\|_{B}$ denotes the Bloch seminorm for the classical Bloch space $\B$ and it is clear that $h \in \B$. Since $\|L \circ \phi\| \leq 1$, take $x_0=0$ in Proposition \ref{propo} and we conclude that $\|(f \circ \phi)'(0)\| \leq \|h\|_B$ for any $\phi \in Aut(B_E)$ so $\|f\|_{inv} < \infty$. \qed \bigskip

\end{document}